\documentclass[accepted]{article}
\usepackage{microtype}
\usepackage{graphicx}
\usepackage{subfigure}
\usepackage{booktabs} 

\usepackage{rotating} 
\usepackage[draft]{hyperref}

\usepackage{icml2018}

\usepackage{lipsum} 

\usepackage{amsmath,amssymb}
\newcommand{\R}{\mathbb{R}} 

\icmltitlerunning{Robust Learning of Trimmed Estimators via Manifold Sampling}

\begin{document}

\twocolumn[
\icmltitle{Robust Learning of Trimmed Estimators via Manifold Sampling}




\begin{icmlauthorlist}
\icmlauthor{Matt Menickelly}{anl}
\icmlauthor{Stefan M.~Wild}{anl}
\end{icmlauthorlist}

\icmlaffiliation{anl}{Mathematics and Computer Science Division, Argonne 
National Laboratory, Lemont, Illinois, USA}
\icmlcorrespondingauthor{Matt Menickelly}{mmenickelly@anl.gov}
\icmlcorrespondingauthor{Stefan M.~Wild}{wild@mcs.anl.gov}

\icmlkeywords{Machine Learning, ICML}

\vskip 0.3in
]

\printAffiliationsAndNotice{}  

\begin{abstract}
We adapt a manifold sampling algorithm for the nonsmooth, nonconvex 
formulations of learning that arise when imposing robustness to outliers present in the 
training data. We demonstrate the approach on objectives based on 
trimmed loss. Empirical results show that the method has favorable scaling 
properties. Although savings in time come at the expense of not 
certifying optimality, the algorithm consistently returns high-quality 
solutions on the trimmed linear regression and multiclass classification 
problems tested.
\end{abstract}

\section{Introduction}
\label{sec:intro}

A frequent challenge in supervised learning is to perform training that is 
robust to outliers in a labeled training dataset 
$\mathbb{T}=\{(x^i,y^i)\}_{i=1}^N\subset 
\mathbb{R}^d \times \mathbb{R}$. 
Here, we address this challenge by using 
\emph{trimmed estimators}. Formally, for 
$l_i(F(x^i,w),y^i)$,
a scalar-valued loss function of the label $y^i$ and model prediction 
$F(x^i,w)$, we seek weights $w\in \mathbb{R}^n$ that solve
\begin{equation}
\min_w \displaystyle\frac{1}{q}\displaystyle\sum_{j=1}^q 
l_{(j)}\left(F(x^{(j)},w),y^{(j)}\right),
 \label{eq:prob2}
\end{equation}
where 
$(q)$ 
denotes
the index associated with the $q$th-order statistic 
(i.e., $l_{(j-1)}\left(F(x^{(j-1)},w),y^{(j-1)}\right)\leq 
l_{(j)}\left(F(x^{(j)},w),y^{(j)}\right)$ for $j=2,\ldots, N$.
Robustness properties of trimmed estimators have been explored in, for example, 
\cite{Rousseeuw1984,Rousseeuw1987}. 
Such estimators have a \emph{high breakdown point}; 
as a particular example in the context of least squares regression,
when $q = \lfloor\frac{N}{2}\rfloor$, the estimator obtained by
solving \eqref{eq:prob2} can resist the adverse effects of up to $50\%$ outliers present in the data. 
Despite such promising properties, however, a barrier to solving \eqref{eq:prob2} has 
been 
its inherent nonconvexity and nonsmoothness (even if each $l_i(\cdot)$ is smooth 
and convex), rendering global solutions intractable for large values of $N$ or 
$d$. A majority of approaches to solving \eqref{eq:prob2} in the 
context of least-squares regression (e.g., \cite{RousseeuwDriessen})
have focused on heuristics that sample subsets of size $q$ from the dataset and 
then follow a sample average approximation-like approach of solving
\eqref{eq:prob2} restricted to this subset.
Such an approach necessitates heuristics because 
complete enumeration would involve the prohibitive evaluation of ${N \choose q}$
many subsets in order to be assured of a set of smooth (possibly convex) problems. 
In the context of linear regression (i.e., $F(x,w) = w^\top x$), recent work
\cite{Bertsimas2014} has
 proposed solving the least quantile regression problem (of 
the form \eqref{eq:prob2} but with the sum starting at $j=q$) to global 
optimality through integer reformulations. 
Although that approach necessarily returns certificates of global optimality, 
integer optimization is not scalable, and least
quantile regression is strictly more statistically inefficient than solving 
\eqref{eq:prob2}; see \cite{Rousseeuw1984}.
Recently, a stochastic proximal gradient method 
targeting general problems of the form
\eqref{eq:prob2} by solving a relaxation of the problem in $n+N$ dimensions
has been proposed \cite{AravkinDavis}.


In this paper we use the observation that the objective in \eqref{eq:prob2}
can be expressed as the composition of a piecewise linear function 
with the mapping $l(w) 
= [\ldots, l_i(F(x^i,w),y^i), \ldots]^\top$ 
in order to apply a \emph{manifold sampling algorithm}. 
Let us make this piecewise linear structure explicit.
Given a training dataset $\mathbb{T}$ and weights $w$, 
the values of $\{l_i(F(x^i,w),y^i)\}_{i=1}^N$ are fixed,
 and the objective in
\eqref{eq:prob2} 
%
is defined by the set of $q$-tuples of active indices 
\begin{equation*}
\begin{array}{rl}
\hspace{-4pt} I^{q,N}(l(w)) =& \Big\{
(i_1,\ldots,i_q): \, 
 l_{i_j}\left(F(x^{i_j},w),y^{i_j}\right) \\
&  \leq l_{(q)}\left(F(x^{(q)},w),y^{(q)}\right)\Big\}.
 \label{eq:active2}
 \end{array}
\end{equation*}
Then, the objective in \eqref{eq:prob2} is expressible as $(h^{q,N} \circ 
l):\mathbb{R}^n\to\mathbb{R}$, 
where $h^{q,N}$ is the \emph{continuous selection} \cite{Scholtes2012}
\begin{equation*}
h^{q,N}\left( l(w) \right) = \left \{g^a\left(l(w)\right) : a\in I^{q,N}(l(w)) 
\right\}
 \label{eq:ref}
\end{equation*}
and the \emph{selection function} $g^a$ for $a\in I^{q,N}(l(w))$ is defined 
componentwise by
\begin{equation}
\label{eq:ha2}
g^a_i\left(l(w)\right) = \left\{
\begin{array}{ll}
 1/q & i\in a \\
 0 & \text{ otherwise. }
\end{array}
\right.
\end{equation}

Even though 
$I^{q,N}(l(w))$ need not be 
singleton, $h^{q,N}(l(w))$ is singleton for all $w$.
%


We observe that
\eqref{eq:prob2} 
is generally a nonconvex, nonsmooth 
optimization problem. 
In particular, even if the loss function $l$ is smooth, 
the objective is potentially nondifferentiable at all values of $w$ where 
 $I^{q,N}(l(w))$ is nonsingleton. 
Similarly, even if the loss function $l$ is convex, 
$h^{q,N}$
is generally a nonconvex piecewise linear function, 
and so \eqref{eq:prob2} falls outside the 
scope of convex composite 
optimization; see \cite{DuchiRuan} and the references therein. 


\section{Manifold Sampling}
\label{sec:MS}
Manifold sampling is an iterative method for minimizing compositions of a 
continuous selection $h^{q,N}:\mathbb{R}^N\to\mathbb{R}$ 
and a smooth function $l:\mathbb{R}^n\to\mathbb{R}^N$.
The case where $h^{q,N}$ is piecewise linear, considered here, is 
analyzed in \cite{KLW17}. 
The manifolds in the case of \eqref{eq:prob2} can be thought of as the 
regions of $\R^d$ within which the objective is piecewise smooth. These 
regions can be catalogued by the set of 
\emph{active indices} 
of $h^{q,N}(l(\cdot))$ (i.e., the set of indices $I^{q,N}(l(\cdot))$).  
At iteration $k$ of a manifold sampling method, there are a current iterate $w^k$
and a trust-region radius $\Delta_k$.
The iteration involves a direction-finding routine and
a \emph{manifold sampling loop}, which we now overview; see
\cite{KLW17,ROBOBOA} for details. 
By the end of the loop, a \emph{generator set} $\mathbb{A}^k$ is determined, at 
least containing $I^{q,N}(l(w^k))$,
and at most containing $\{I^{q,N}(l(y)) : y\in\mathcal{B}(w^k,\Delta_k)\}$.
Following the strategy in \cite{ROBOBOA},
the step $d^k$ is taken from the approximate solution to 
\begin{equation*}
\small
\hspace{-4pt}
\begin{array}{l}
\displaystyle\min_{\tau\in\mathbb{R},d\in\mathbb{R}^n} \tau: \|d\|_2 \leq \Delta_k, \\
g^a(l(w^k)) - h^{q,N}(l(w^k)) + 
\nabla(g^a \circ l)(w^k)^\top d \leq \tau \quad \forall a\in \mathbb{A}^k,
\end{array}
\end{equation*}
where $g^a$ is the selection function (see \eqref{eq:ha2}) corresponding to 
$a\in \mathbb{A}^k$.
The iteration then proceeds similarly to a standard trust-region method.

We propose a stochastic variant by making two adjustments to the manifold 
sampling algorithm. First, given a sample $S\subseteq\{1,\dots,N\}$, we 
consider a \emph{sampled} smooth function
$l_S(w): \mathbb{R}^n\to\mathbb{R}^{|S|}$;
intuitively, $l_S(w)$ is the projection of $l(w)$ onto the coordinates 
corresponding to indices in $S$. 
Then, in the $k$th iteration, the output of the manifold sampling loop can be 
regarded as a function of a sample $S^k$,
a continuous selection $h$ (in particular, $h^{\lfloor (q/N)|S^k| \rfloor,|S^k|}$),
 the current iterate $w^k$, and the trust-region radius $\Delta_k$ (i.e., we 
write $(\tau_k,d^k) = \mathcal{M}(S^k,h,w^k,\Delta_k)$). 
Second, we modify the step acceptance test 
so that given a secondary sample $S^{k'}$, 
continuous selection $h$ (here, $h^{\lfloor (q/N)|S^{k'}| \rfloor,|S^{k'}|}$), 
and current iterate $w^k$, 
\begin{equation*}
\small
\rho_k(S^{k'},h,w^k,\tau_k,d^k) = \displaystyle\frac{h(l_{S^{k'}}(w^k))- 
h(l_{S^{k'}}(w^k+d^k))}{-\tau_k}.
\label{eq:rho}
\end{equation*}
A statement of the algorithm is given in Algorithm~\ref{alg:sms}. 

\begin{algorithm}[tb]
   \caption{(Sample-Based) Manifold Sampling}
   \label{alg:sms}
\begin{algorithmic}
   \STATE {\bfseries Input:} 
   parameters $\gamma_{\rm dec}\in(0,1)$,
   $\gamma_{\rm inc}>1$, $\eta\in(0,1)$, initial point $w^0\in\mathbb{R}^n$, 
    trust-region radius $\Delta_0>0$, and counter $k=0$
   \REPEAT
   \STATE Sample $S^k\subset \{1,\dots,N\}$
   \STATE $(\tau_k,d^k) \gets \mathcal{M}(S^k,h^{\lfloor q|S^k| 
\rfloor,|S^k|},w^k,\Delta_k)$
   \STATE Sample $S^{k'}\subset\{1,\dots,N\}$
    \IF{$\rho_k(S^{k'},h^{\lfloor q|S^{k'}| 
\rfloor,|S^{k'}|},w^k,\tau_k,d^k) > \eta$}
   \STATE $w^{k+1}\gets w^k + d^k$
   \STATE $\Delta_{k+1}\gets\gamma_{\rm inc}\Delta_k$
   \ELSE
   \STATE $w^{k+1}\gets w^k$
   \STATE $\Delta_{k+1}\gets\gamma_{\rm dec}\Delta_k$
   \ENDIF
   \STATE{$k\gets k+1$}
   \UNTIL{budget exhausted}
\end{algorithmic}
\end{algorithm}

Although some manifold sampling variants do not require gradients $\nabla l$, 
these gradients are readily available in the problems tested here, and hence 
we make direct use of them (but not second-order information $\nabla^2 l$). 
Consequently, we remark that because we are not using model/approximate 
Hessians, 
there is no need for the ``acceptability test'' (relative to a stationary 
measure) discussed in \cite{KLW17,ROBOBOA}.
We also remark that while the convergence analysis in \cite{KLW17,ROBOBOA}
applies to the deterministic case wherein the sample $S^k=\{1,\dots,N\}$ on every iteration $k$, 
a proper analysis of the stochastic variant where $|S^k|<N$ is beyond the scope  
of this paper. 
In future work, we will seek to establish convergence of the stochastic variant 
of Algorithm~\ref{alg:sms} to Clarke stationary points with probability
1 (given an appropriate sampling rate). 

\section{Numerical Results}
\label{sec:num_results}

We now demonstrate the proposed approach on two standard problems.
We implemented Algorithm~\ref{alg:sms} in \texttt{MATLAB}. 
Throughout, we initialize $w^0=0$ and $\Delta_0=10$
and choose $\gamma_{\rm inc}=1.01, \gamma_{\rm dec} = 0.99, \eta = 10^{-3}$. 

\begin{table*}[th]
\caption{Average (over 30 trials) results of ability to correctly 
identify outliers on the regression experiments.\label{table:regression}}
\begin{center}
\begin{small}
\begin{footnotesize}
\begin{tabular}{l|cccc|cccc|cccc|cccc}
%
%
%
$d$ & 5 & 5 & 5 & 5 & 10 & 10 & 10 & 10 & 20 & 20 & 20 & 20 & 100 & 100 & 100 & 
100 \\
$\frac{p}{1000}$ & 2 & 5 & 10 & 50 & 2 & 5 & 10 & 50 & 2 & 5 & 10 & 50 & 2 & 5 
& 10 & 50
\\ \midrule \midrule
\multicolumn{1}{c}{} & \multicolumn{16}{c}{\textbf{LQS}} \\
TPR & 
100 & 100 & 100 & 100 & 100 & 100 & 100 & 100 & 65.6 & 67.1 & 66.8 & 66.0 & 
43.8 & 43.5 & 43.9 & 43.9\\
 FPR & 0 & 0 & 0 & 0 & 0 & 0 & 0 & 0 & 23.0 & 22.0 & 22.1 & 22.7 & 37.6 & 37.7 
& 37.9 & 37.9\\
 time & 0.1 & 0.3 & 0.6 & 2.6 & 0.2 & 0.4 & 0.8 & 4.0 & 0.2 & 0.6 & 1.0 & 4.9 & 
2.4 & 4.4 & 9.0 & 42.0\\
 \midrule
\multicolumn{1}{c}{} & \multicolumn{16}{c}{\textbf{DMS}} \\
 TPR & 100 & 100 &  100 &  100 &  100 &  100 &  100 &  100 &  100 &  100 &  100 
&  100 &  100 &  100 &  100 &  100 \\
 FPR & 0 & 0 & 0 & 0 & 0 & 0 & 0 & 0 & 0 & 0 & 0 & 0 & 0 & 0 & 0 & 0\\
 time & 8.5 & 9.9 & 9.9 & 31.2 & 10.6 & 10.1 & 12.6 & 38.2 & 15.1 & 13.7 & 18.4 
& 45.7 & 37.9 & 47.4 & 91.9 & 436\\
\midrule
\multicolumn{1}{c}{} & \multicolumn{16}{c}{\textbf{SMS}} \\
 TPR & 100 & 100 &  100 &  100 &  100 &  100 &  100 &  100 &  100 &  100 &  100 
&  100 &  100 &  100 &  100 &  100 \\
 FPR & 0 & 0 & 0 & 0 & 0 & 0 & 0 & 0 & 0 & 0 & 0 & 0 & 0 & 0 & 0 & 0\\
 time & 0.9 & 1.0 & 1.1 & 2.2 & 1.2 & 1.4 & 1.6 & 2.9 & 1.7 & 2.0 & 2.3
& 4.6 & 2.5 & 4.6 & 5.3 & 15.0\\
\bottomrule
\end{tabular}
\end{footnotesize}
\end{small}
\end{center}
\end{table*}

\paragraph{Regression Problems.}
For the first set of experiments, we follow 
\cite{RousseeuwDriessen,Bertsimas2014} and consider a trimmed form of linear 
regression, where the loss function is absolute loss:
\begin{equation}
\min_w \displaystyle\frac{1}{q}\displaystyle\sum_{i=1}^q |w^\top x^{(i)} - 
y^{(i)}|.
 \label{eq:abs_lin_reg}
\end{equation} 
Although we focused on smooth functions 
$l(w)$ for ease of presentation, it is straightforward to  
consider the manifolds associated with both the quantiles and those induced
by the absolute value operator; see \cite{LMW16}. 

We randomly generate synthetic data and then contaminate a random subset of 
that data to simulate outliers, using the same 
method as in \cite{RousseeuwDriessen,Bertsimas2014}.
In particular, given a dimension $d$ and a number of samples $N$, we independently generate 
$\{x^i_j\}_{j=1,\dots,d}^{i=1,\dots,N} \sim \mathcal{N}(0,100)$. 
Letting $e\in\mathbb{R}^d$ denote a vector of all ones, we randomly generate
$y^i = e^\top x^i + \sigma^i$, where $\sigma^i\sim \mathcal{N}(0,1)$. 
We then corrupt $40\%$ of the dataset by both
(i) randomly selecting $\lfloor 0.2N\rfloor$ indices $i$ and replacing the first
entry $x^i_1$ with a random entry drawn from 
$\mathcal{N}(100,100)$ and (ii) randomly selecting an additional $\lceil 
0.2N \rceil$ indices $i$ and deterministically replacing
$y^i$ with $y^i+1000$. This contamination results in both \emph{bad leverage 
points} and \emph{vertical outliers} in the dataset
\cite{RousseeuwvanZomeren}.
For each experiment, we set $q=\lfloor 0.6N\rfloor$ and generate 
30 sets of data. 

We compare both a deterministic and stochastic implementation of Algorithm~\ref{alg:sms} to the \texttt{lqs} 
function of the popular \texttt{MASS} package \cite{MASS} in \texttt{R}, which 
uses a heuristic to approximately solve \eqref{eq:abs_lin_reg}.
It was demonstrated in \cite{Bertsimas2014} that the \texttt{lqs} function 
tends to return solutions quickly but that the solutions are bad estimates of 
global minima. Here, we highlight a related concern:
for ``larger'' problems (especially in terms of $d$), the solutions returned by 
\texttt{lqs} do not perform well in terms of true positive rate (TPR) and false 
positive rate (FPR), where ``positive'' refers to whether or not a given 
data point is an outlier.

We terminate the deterministic variant (``DMS") of our method (i.e.,  
$S^k=S^{k'}=\{1,\dots,N\}$ in each iteration) based on 
the stopping criterion $\Delta_k < .01$
and we terminate the stochastic variant (``SMS") of our method (where we draw 
random samples without replacement of size 
$|S^k|,|S^{k'}|=\max\{0.01 p, (10^{-6}p)/\Delta_k^4\}$) after 100 epochs 
(that is, we budget $100p$ draws from the data).
Both tested variants consistently achieve a TPR of 100\% and a FPR of 0\%.
Consequently, in these experiments our method effectively returned the global 
minimum of \eqref{eq:abs_lin_reg}: having identified the 
outliers correctly, the solution of a single linear program in postprocessing 
would return a global solution $w^*$.

Table~\ref{table:regression} shows, for a variety of $d$ and $p$, the 
TPR, FPR, and the time (in seconds) to return a solution for \texttt{lqs}, DMS, and SMS.
We note that we examine larger datasets than the largest 
considered in \cite{RousseeuwDriessen} and \cite{Bertsimas2014}, which are 
$(d=5, p=50000)$ and $(d=20,p=10000)$, respectively. 
We also 
note that our times for identifying what are tantamount to global minima 
are significantly less than what is reported 
for the integer formulations for least quantile regression considered in 
\cite{Bertsimas2014}; 
the disadvantage of our approach is that our method cannot provide a 
certificate of global optimality.
We recall that \cite{Bertsimas2014} focused on the least quantile formulation, 
as opposed to the trimmed estimator formulation; from a combinatorial
perspective, trying to solve \eqref{eq:abs_lin_reg} to global optimality likely 
would further increase solution time.

 \begin{figure*}[ht]
\vskip 0.2in
\begin{center}
\centerline{\includegraphics[width=.333\textwidth]{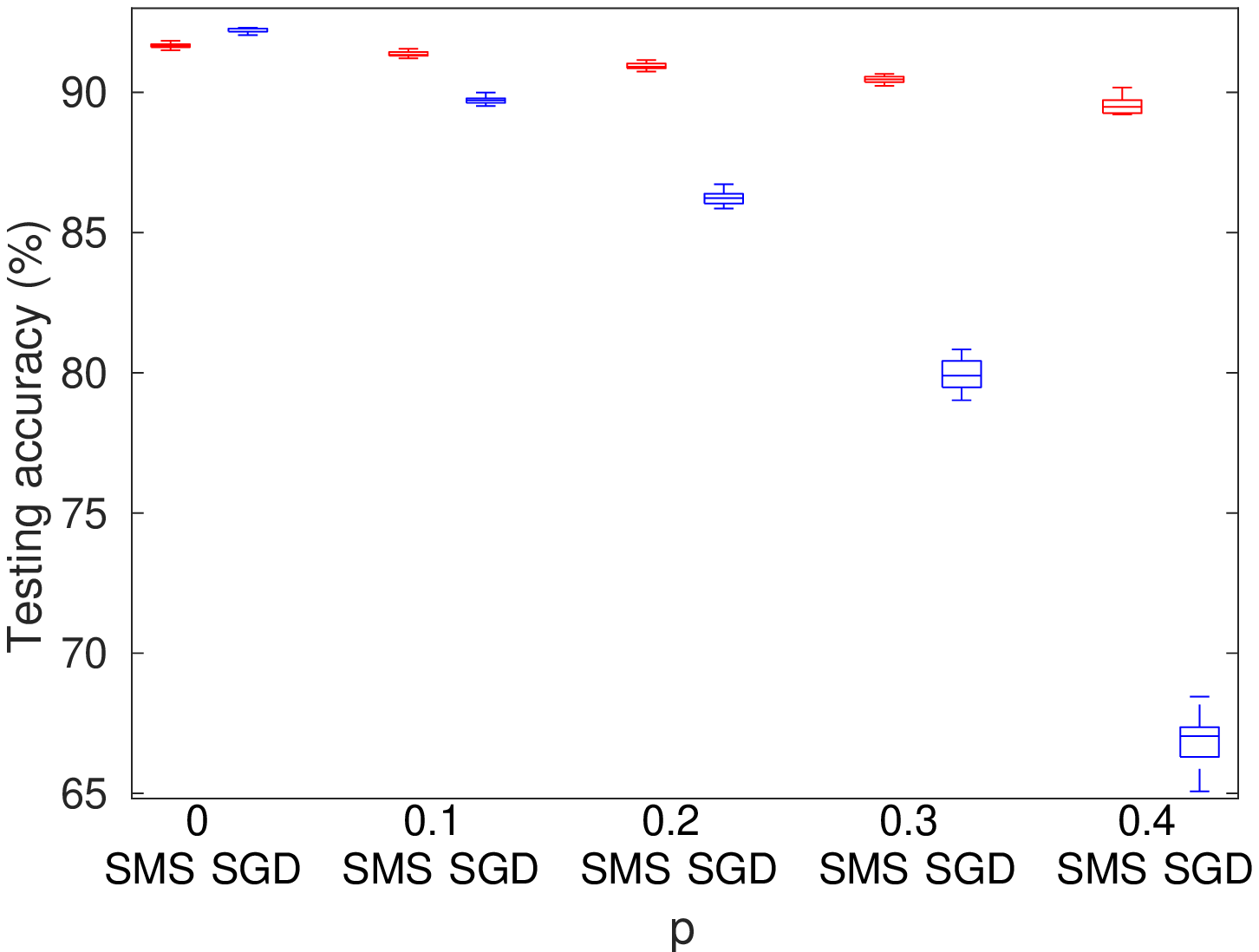}
\includegraphics[width=.333\textwidth]{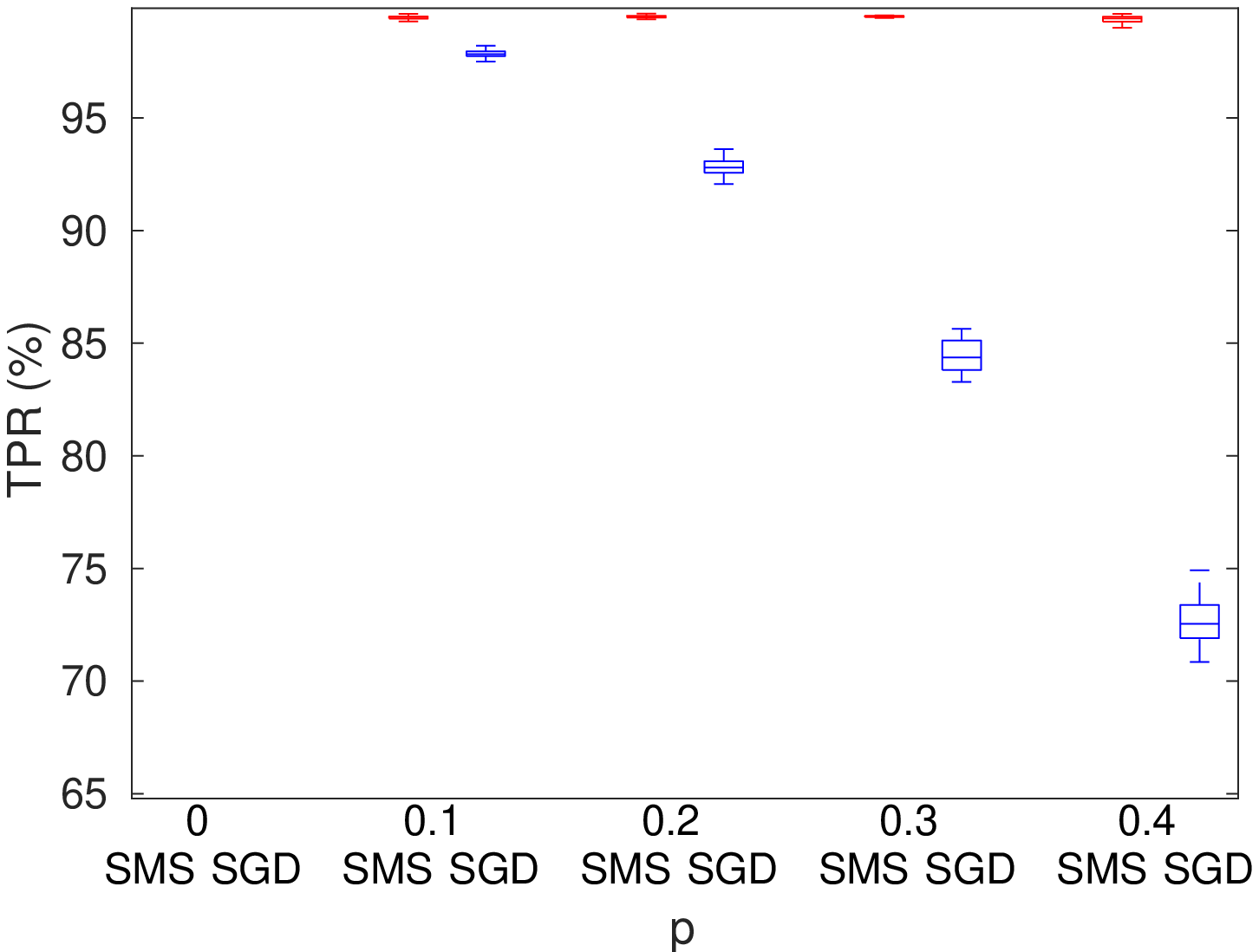}
\includegraphics[width=.333\textwidth]{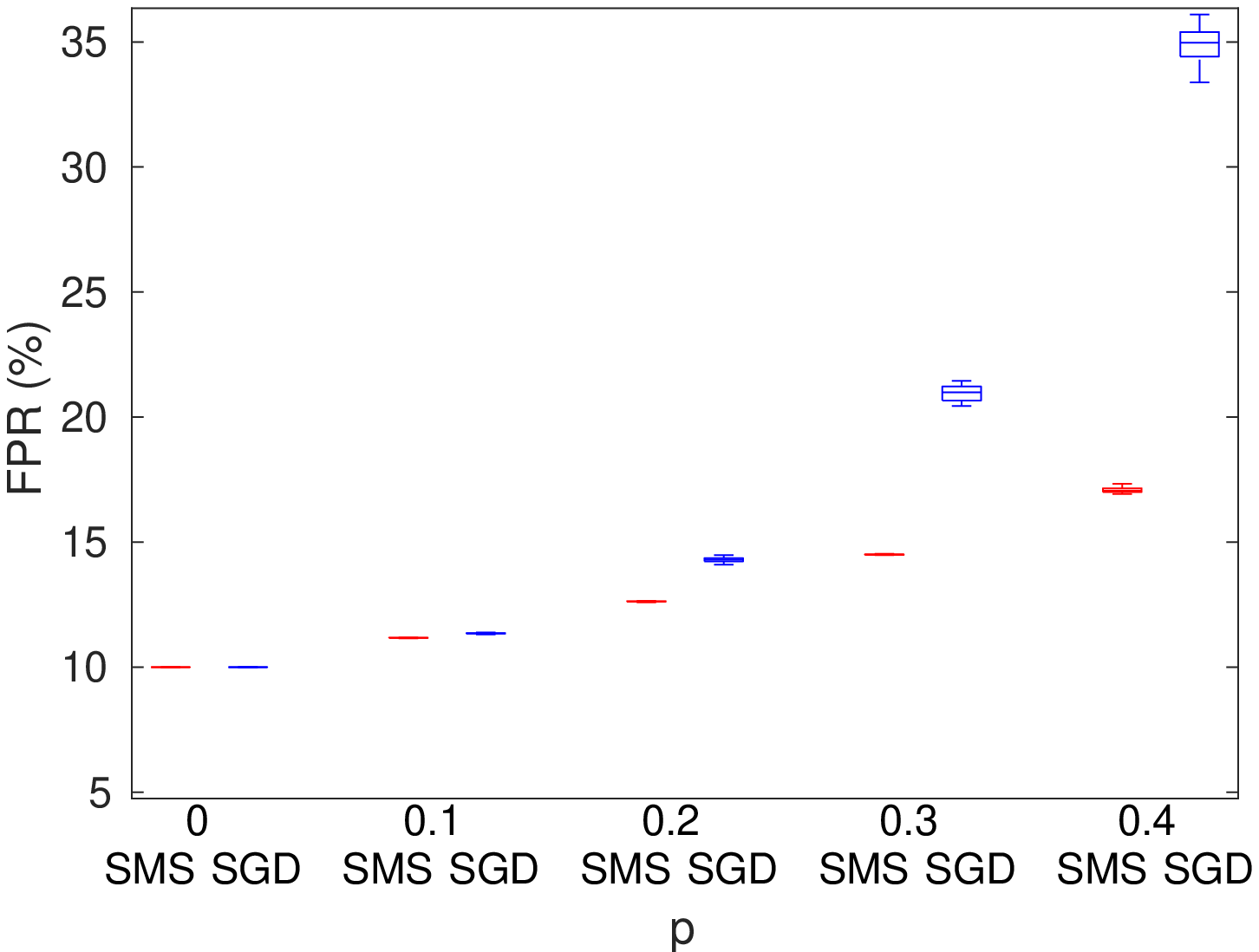}
}
\centerline{\includegraphics[width=.333\textwidth]{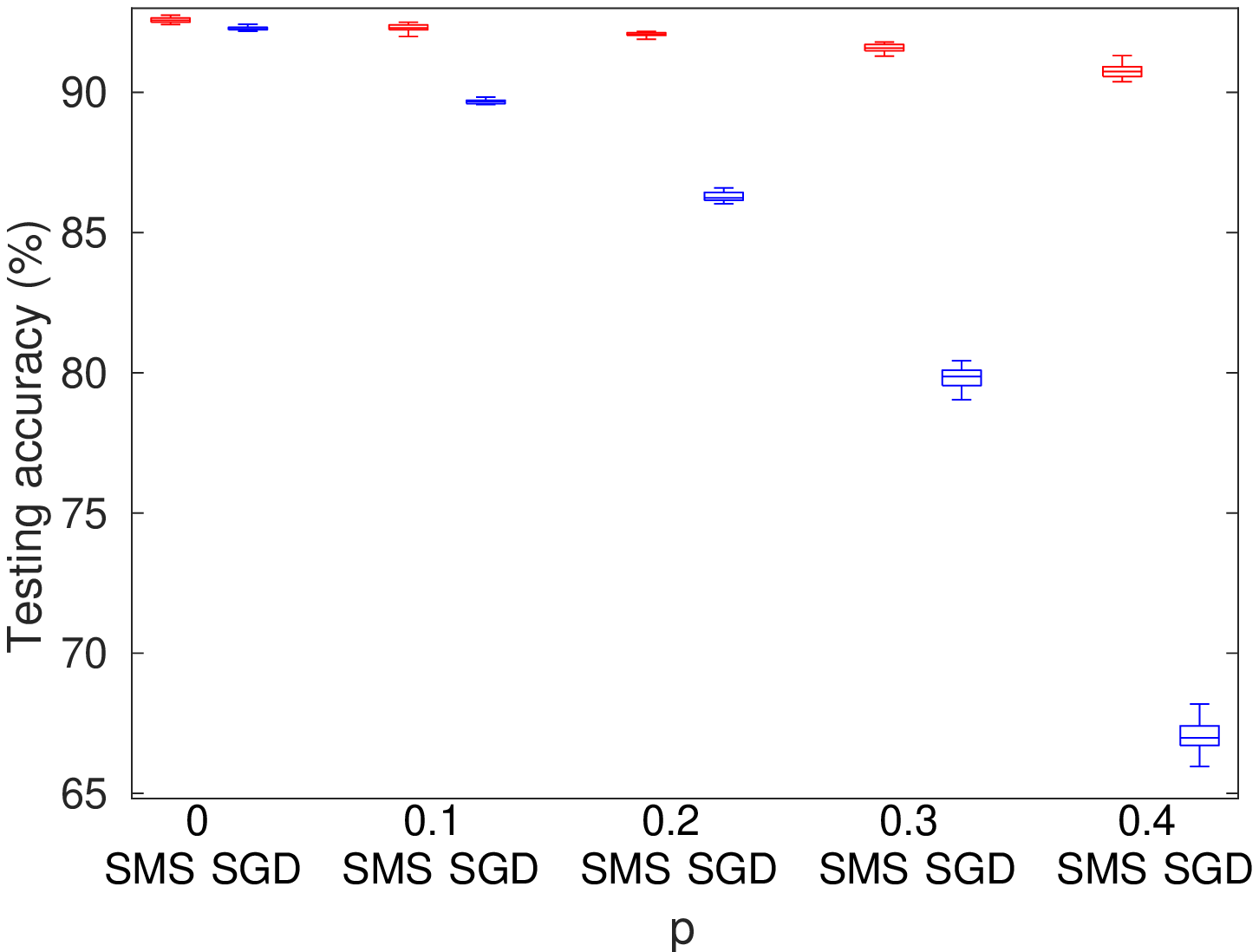}
\includegraphics[width=.333\textwidth]{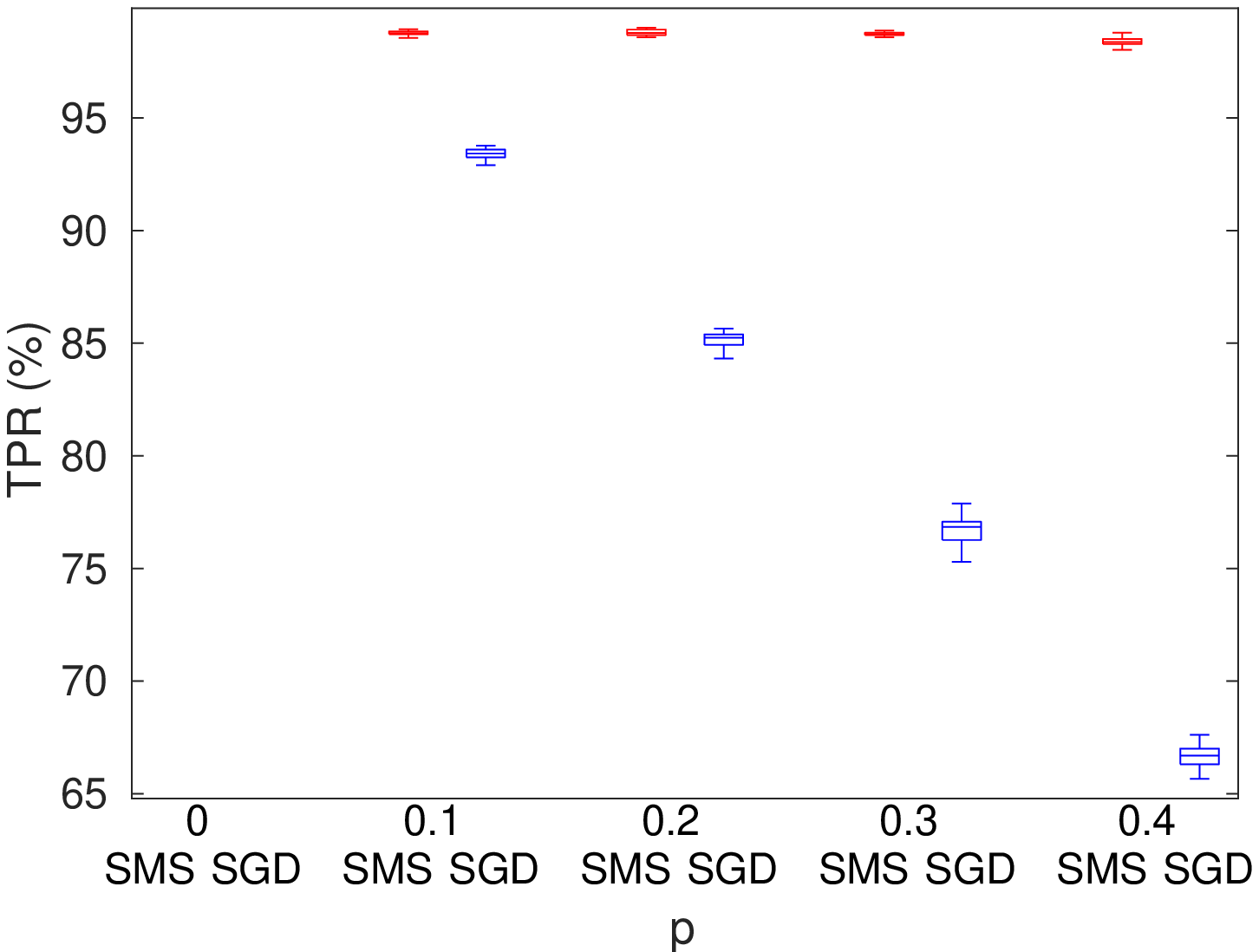}
\includegraphics[width=.333\textwidth]{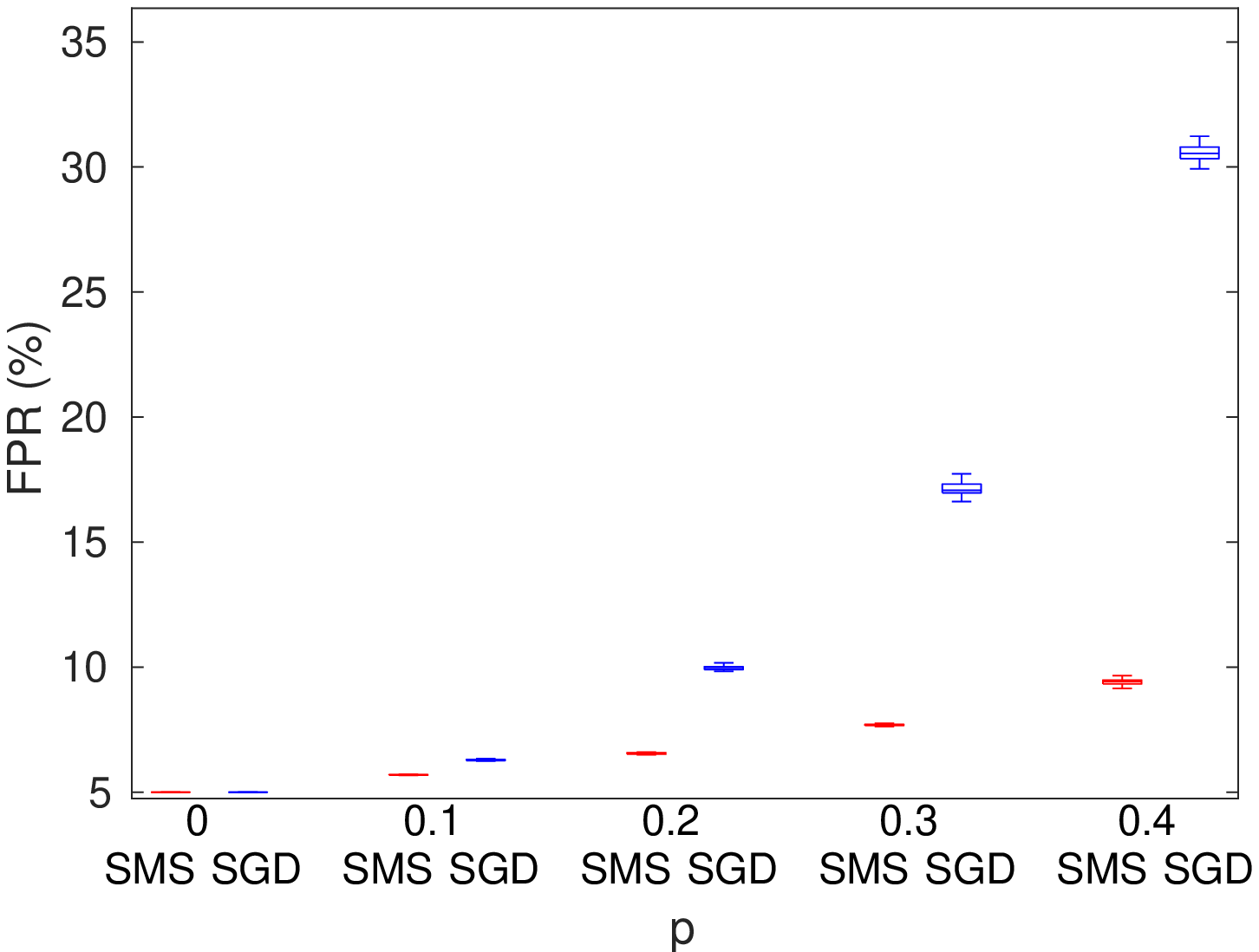}
}
\caption{30 instances of stochastic implementation (SMS) of 
Algorithm~\ref{alg:sms} and standard SGD approach on the MNIST classifier 
problem.
\textbf{Top row:} $q = \lfloor(1-p-0.1)N\rfloor$, \textbf{bottom row:} $q = \lfloor(1-p-0.05)N\rfloor$.
\textbf{Left:} testing accuracy, 
\textbf{center:} TPR on training dataset, 
\textbf{right:} FPR on training dataset.
Because of the experimental setup, the lower bound on FPR of any 
classifier is $10\%$ or $5\%$ for the first and second value of $q$, 
respectively. At the $p=0$ level of contamination, TPR is 
undefined (there are no true outliers).}
\label{fig:mnist_results}
\end{center}
\vskip -0.2in
\end{figure*}

\paragraph{Multiclass Classification on MNIST.}
Following \cite{AravkinDavis}, we apply the stochastic variant of 
Algorithm~\ref{alg:sms} to the well-known task
of classifying digits $0$ through $9$ in the MNIST dataset ($N=60000$) \cite{MNIST}. 
Each data point 
$(x^i,y^i) \in \mathbb{R}^{785}\times\mathbb{R}^{10}$, where $y^i$ is a 
``one-hot'' vector indicating the class membership. 
The learning is performed with the multiclass softmax function
$\ell_i:\mathbb{R}^{10\times 785}\to\mathbb{R}$ defined by
 $$\ell_i(W)=-\log\left(\displaystyle\frac{(y^i)^\top\exp(W 
x^i)}{e^\top\exp(W x^i)}\right),$$
where $\exp$ is taken entrywise and $e$ is the vector of ones. 
Class predictions for any input $x\in\mathbb{R}^{785}$ are then made by evaluating 
the product $Wx\in\mathbb{R}^{10}$ and predicting the class as that corresponding to the 
index of the maximum entry.
We consider the trimmed estimator resulting from 
\eqref{eq:prob2}, with $w$ denoting the vector form of $W$.
 To artificially produce outliers, we contaminate the labels of the data by 
choosing a fraction
 $p\in[0,1]$ and changing the label $y^i$ on a randomly selected subset of size $\lceil pN \rceil$ 
 of the training data to $(y^i + 1)\mod 10$.
 
 We compare the stochastic variant of Algorithm~\ref{alg:sms} (``SMS'') with a 
stochastic gradient descent (SGD) method, which has well-tuned step-size 
and batch-size parameters but 
 solves the untrimmed formulation (i.e., \eqref{eq:prob2} with 
$q=N=60000$). 
 For various levels of $p$, when applying SMS, we consider 
 an overestimation of the number of outliers in the data (i.e., we select either
 $q = \lfloor (1-p-0.1)N \rfloor$
 or $q = \lfloor (1-p-0.05)N \rfloor$). This simulates a 
 realistic application of an outlier detection method, since the proportion of outliers $p$ generally
 cannot be known exactly a priori but can, one hopes, be upper bounded. 
 For SMS, we used a sample on each iteration of size
 $|S^k|=|S^{k'}|=\max\{40,(10^{-6}/\Delta_k^4)\}$. 
 We stop both methods after 40 epochs (i.e., $40N$ draws are made from the 
training data).
 
Figure~\ref{fig:mnist_results} displays the testing accuracy -
 the accuracy of predictions made by learned classifiers $W^*$ on the uncontaminated
testing MNIST set of size $N=10000$ - and the TPR and FPR
of $W^*$ applied to the contaminated training dataset. 
Recall that ``positive'' here refers to a data point being an outlier;
for the output $W^*$ returned by either method, we label all data points 
corresponding to $l_{(q+1)}(W^*)\leq\dots\leq l_{(N)}(W^*)$
as outliers in postprocessing. 
Since both methods are stochastic, we performed 30 trials of the experiment 
with different random seeds. 
 Because these problems are nonconvex and likely have many local minima, 
 we show the 25th, 50th, and 75th percentiles of the results in boxes and the 
outliers in whiskers. 
The results illustrate how, even with the agnostic initialization $w^0=0$, SMS 
has a tendency to identify a robust solution in terms of the natural 
performance metrics.
Comparisons with a well-tuned stochastic (sub)gradient method may be found in 
the Supplementary Material.

The median results of SMS are comparable to the results of the SMART algorithm reported in 
\cite{AravkinDavis}, 
 but we remark that SMS works in $\mathbb{R}^{7850}$,
 while the SMART approach works in $\mathbb{R}^{7850 + 60000}$. 

\clearpage
\section*{Acknowledgements}
This material was based upon work supported by the applied mathematics program 
of the U.S.\ Department of Energy, Office of Science, Office of Advanced 
Scientific Computing Research,  under Contract DE-AC02-06CH11357.
%

\bibliographystyle{icml2018}

\clearpage
\appendix

\section{Supplementary Material}
It is of interest to compare our method to
a simple \emph{stochastic subgradient descent (SSGD)} method.
Given $w^k\in\mathbb{R}^{n}$, 
we can compute a particular subgradient of \eqref{eq:prob2}
by choosing a fixed element of $I^{\lfloor q|S^k|\rfloor,|S^k|}(l(w^k))$
by using a lexicographic ordering; we label this element (which is a 
$q$-tuple) $I^k$. 
We remark that in practice, and because of the almost-everywhere 
differentiability of $h^{\lfloor q|S^k|\rfloor,|S^k|}\circ l(\cdot)$,
 it is often the case that $I^{\lfloor q|S^k|\rfloor,|S^k|}(l(w^k))$ consists 
of a single $q$-tuple; in this case $I^k$ is well defined with no need for a 
lexicographic ordering. 

The $q$-tuple $I^k$ defines a particular subgradient of $h^{\lfloor 
q|S^k|\rfloor,|S^k|}\circ l$ at $w^k$, namely,
\begin{equation}
\label{eq:subgrad}
g_{S^k}(w^k) \triangleq \frac{1}{\lfloor q|S^k|\rfloor} \displaystyle\sum_{i_j \in I^k} \nabla l_{i_j} \left( F(x^{i_j},w^k),y^{i_j}\right).
\end{equation}
With a (sub)gradient defined, we can employ the SSGD method  
described in Algorithm~\ref{alg:ssgd}.

\begin{algorithm}[b]
   \caption{Stochastic Subgradient Descent (SSGD) Method}
   \label{alg:ssgd}
\begin{algorithmic}
   \STATE {\bfseries Input:} 
   algorithmic parameters 
   initial point $w^0\in\mathbb{R}^{n}$,
   initial step size $\alpha>0$, 
   decay parameter $\delta \geq 0$, 
   and
   batch-size parameter $s \geq 1$.
  \STATE $k \gets 0$
   \REPEAT
   \STATE Randomly sample $S^k\subset \{1,\dots,N\}$ so that $|S^k|=s$. 
   \STATE Compute $g_{S^k}(w^k)$ as in \eqref{eq:subgrad}.
   \STATE $w^{k+1} \gets w^k - \displaystyle\frac{\alpha}{k^{\delta}}g_{S^k}(w^k)$.
   \UNTIL{budget exhausted}
\end{algorithmic}
\end{algorithm}

A deterministic (i.e., $|S^k|=N$ on every iteration $k$) 
variant of such a method was used in \cite{Bertsimas2014} as a means of 
warmstarting
their integer programming model. 
Under modest assumptions, deterministic methods of the form in 
Algorithm~\ref{alg:ssgd} can be shown to converge to minimizers of convex 
(almost-everywhere differentiable) functions; see, for instance, 
\cite{Shor1985}. 
However, most convergence results for methods of nonconvex nonsmooth minimization 
require some form of ``stabilized descent" direction, like those generated in 
bundle 
or gradient-sampling methods (see \cite{makela2001} and \cite{BCLOS18} for
respective literature reviews). 
Manifold sampling methods more closely resemble these latter two classes of 
methods.

\begin{figure*}
\vskip 0.2in
\begin{center}
\includegraphics[width=.49\textwidth]{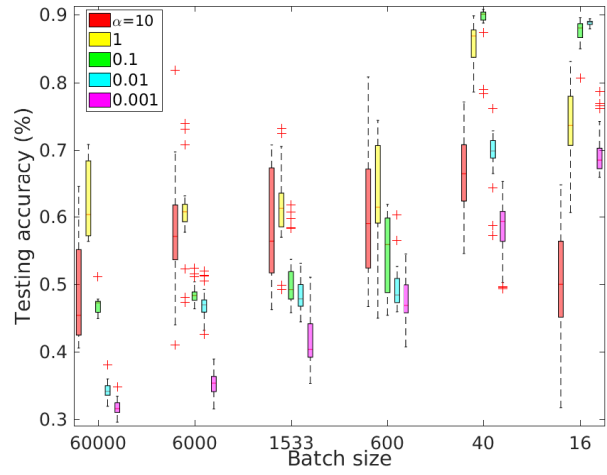}
\caption{Testing accuracy of classifiers learned by SSGD for $p=0.3$ and 
$q=\lfloor (1-p-0.05)N\rfloor$ in the MNIST experiments
over a variety of $(\alpha,s)$ parameter pairs.}
\label{fig:tuning}
\end{center}
\vskip -0.2in
\end{figure*}

In light of this lack of convergence guarantees for Algorithm~\ref{alg:ssgd},
 we performed a tuning of the three parameters, $(\alpha,\delta,s)$,
of Algorithm~\ref{alg:ssgd}.
An illustration of results from our tuning experiments is shown in 
Figure~\ref{fig:tuning}.
Here, we repeat the previously described MNIST experiment with $p=0.3$ and $q=\lfloor (1-p-0.05)N\rfloor$,
 but we apply the SSGD method
to the randomly contaminated problems instead, varying the parameters
$\alpha\in\{10^1,10^0,10^{-1},10^{-2},10^{-3}\}$ and 
$s\in\{N, N/10, N^{2/3}, N/100, N^{1/3}, 
N^{1/4}\}=\{60000, 6000, 1533, 600, 40, 16\}$.
In our experience, setting $\delta>0$ generally deteriorates performance for 
any fixed level of $\alpha, s$.
Hence, in all the results shown, we set $\delta=0$; that is, we effectively 
have a fixed step size.
 
We observe from Figure~\ref{fig:tuning} that parameter tuning
has a significant impact on the empirical performance of the SSGD 
method. Arguably, only three
of the $(\alpha,s)$ pairs (corresponding to $(\alpha,s) \in \{(0.1,40),(0.1,16),(0.01,16)\}$) produce a 
reasonable range of testing accuracy values across the randomly
contaminated instances. Using the best median value as a selection criteria 
suggests that ($\alpha=0.1$, $s=40$) is good for these instances.

\begin{figure*}[t!]
\vskip 0.2in
\includegraphics[width=.49\textwidth]{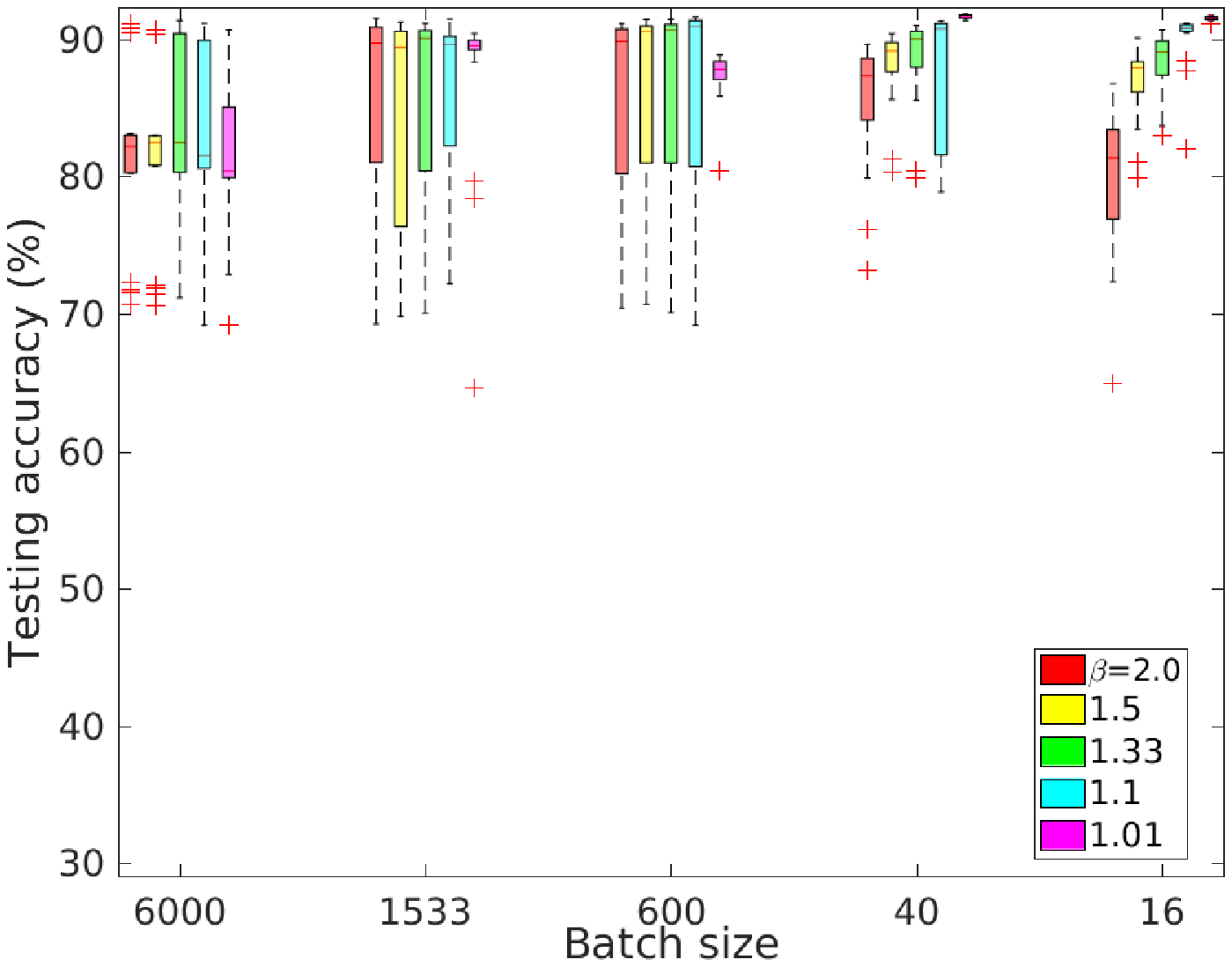}
\includegraphics[width=.49\textwidth]{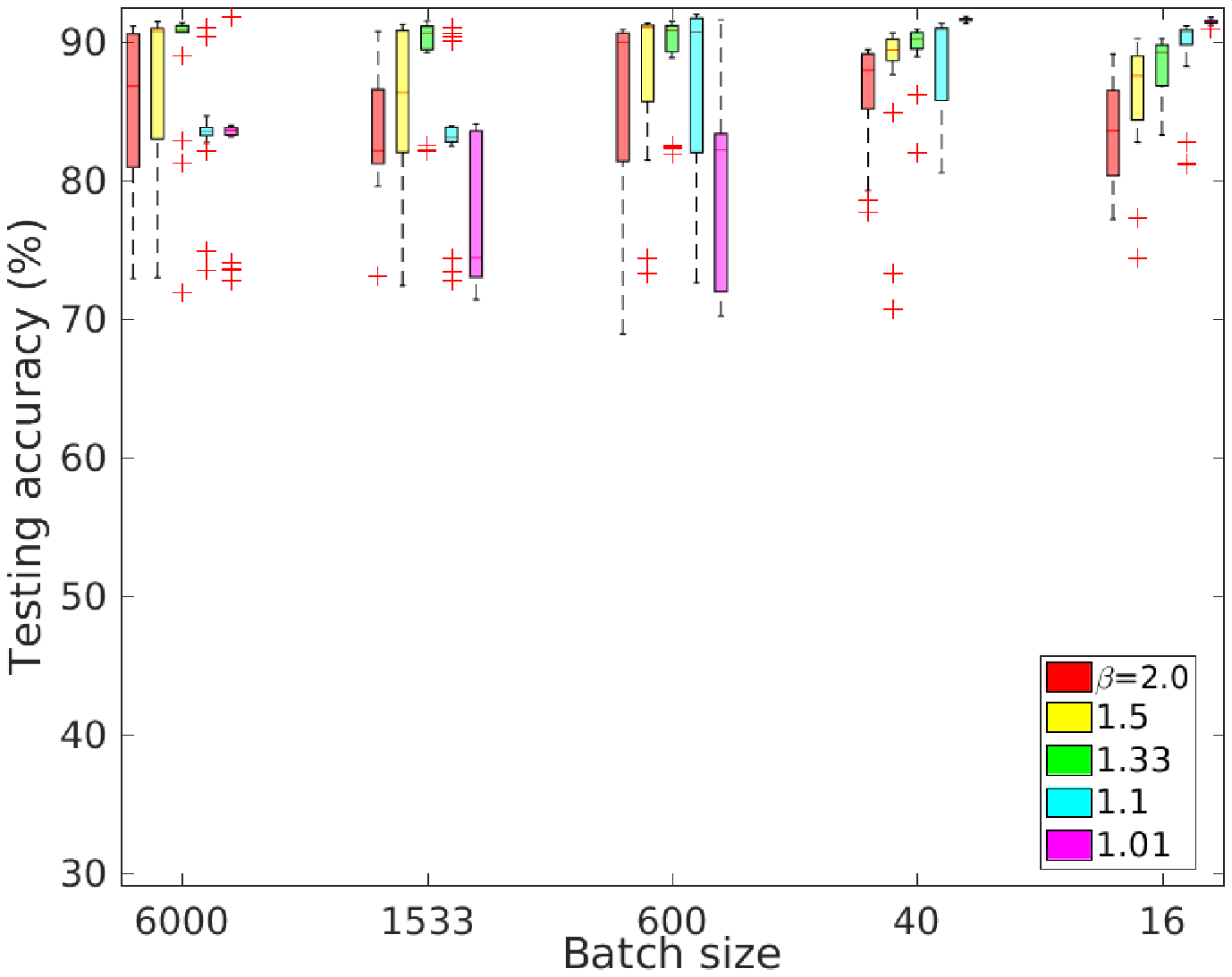}
\caption{Testing accuracy of classifiers learned by SMS for $p=0.3$ and 
$q=\lfloor (1-p-0.05)N\rfloor$ in the MNIST experiments
over a variety of $(\beta,|S^k|=|S^{k'}|)$ parameter pairs. 
\textbf{Left: $\Delta_0=10$, right: $\Delta_0=1$}}
\label{fig:tuning2}
\vskip -0.2in
\end{figure*}

In contrast, 
the empirical performance of SMS is less sensitive to a selection of its parameter 
values. SMS has more parameters to tune, but we observed that the most variance in 
the performance metrics in these MNIST experiments come from the selection of
the ``batch size" $|S^k|=|S^{k'}|$, and the choice of 
the trust-region expansion and contraction parameters,
$\gamma_{\rm inc}$ and $\gamma_{\rm dec}$. 
In Figure~\ref{fig:tuning2} we perform a similar tuning experiment as done for SSGD,
except we tune two parameters $\beta$ and $|S^k|=|S^{k'}|$. For a single value of $\beta$, we define
$\gamma_{\rm inc}=\beta$ and $\gamma_{\rm dec}=1/\beta$. 
We show results for all possible combinations of $\beta\in\{2.0,1.5,1.33,1.1,101\}$,
constant batch sizes in $\{6000,1533,600,40,16\}$, and $\Delta_0\in\{10,1\}$. 
We remark that we did not test a batch size of $60000$ only because
it is typical for SMS to run for only one or two iterations of Algorithm~\ref{alg:sms} in this setting, since we only set a budget of 40 epochs. 

\begin{figure*}[t!]
\begin{center}
\centerline{\includegraphics[width=.333\textwidth]{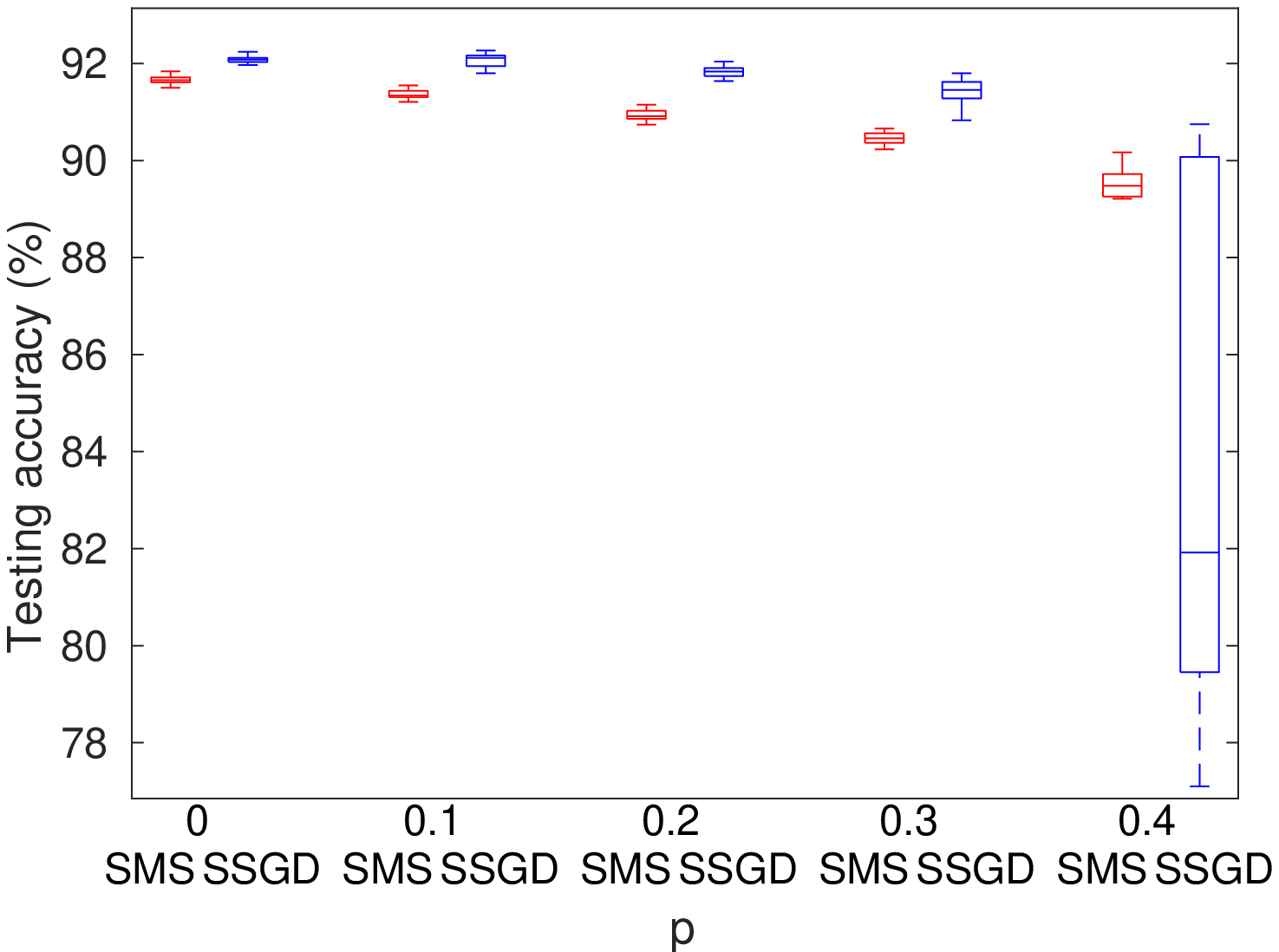}
\includegraphics[width=.333\textwidth]{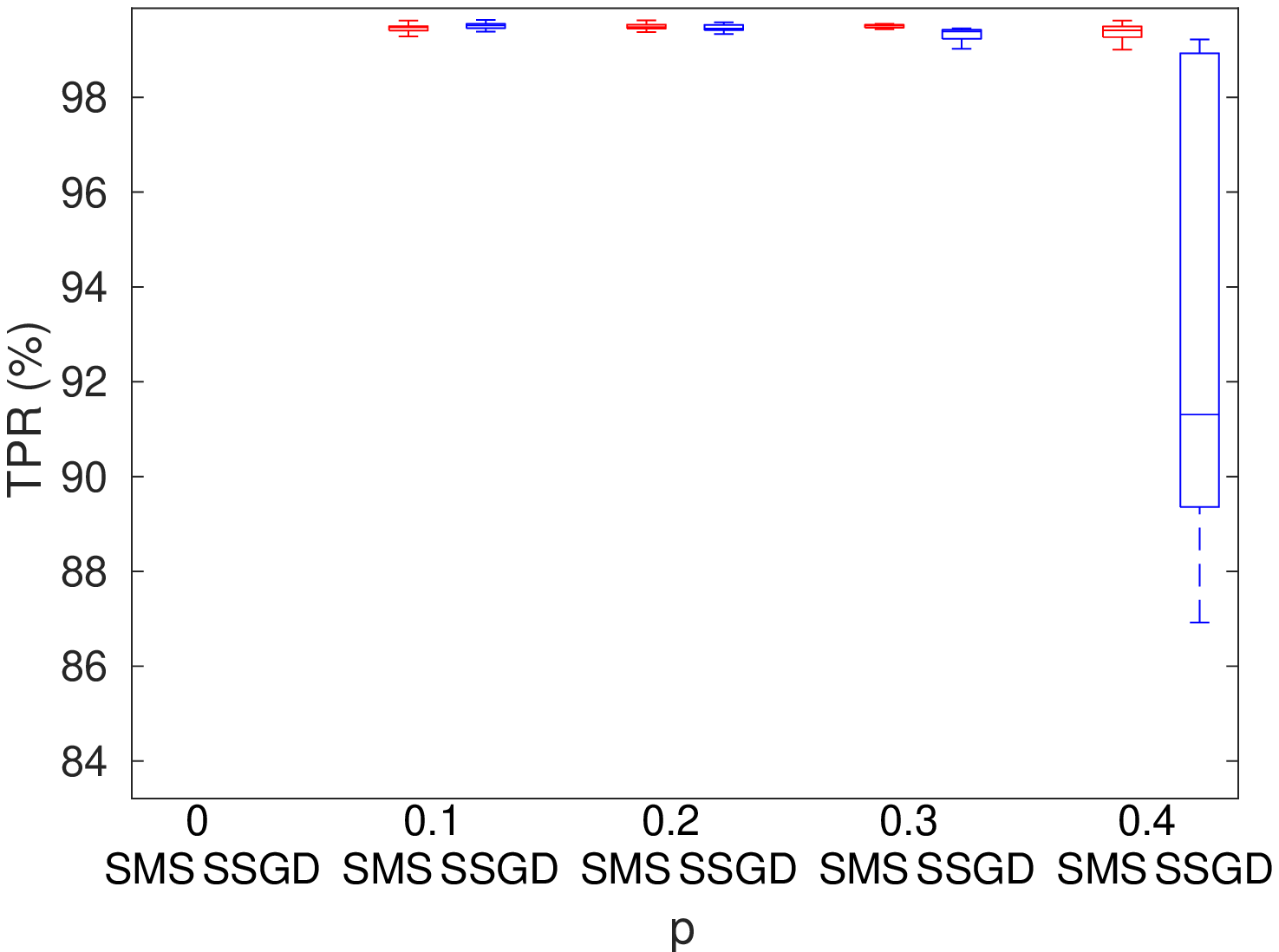}
\includegraphics[width=.333\textwidth]{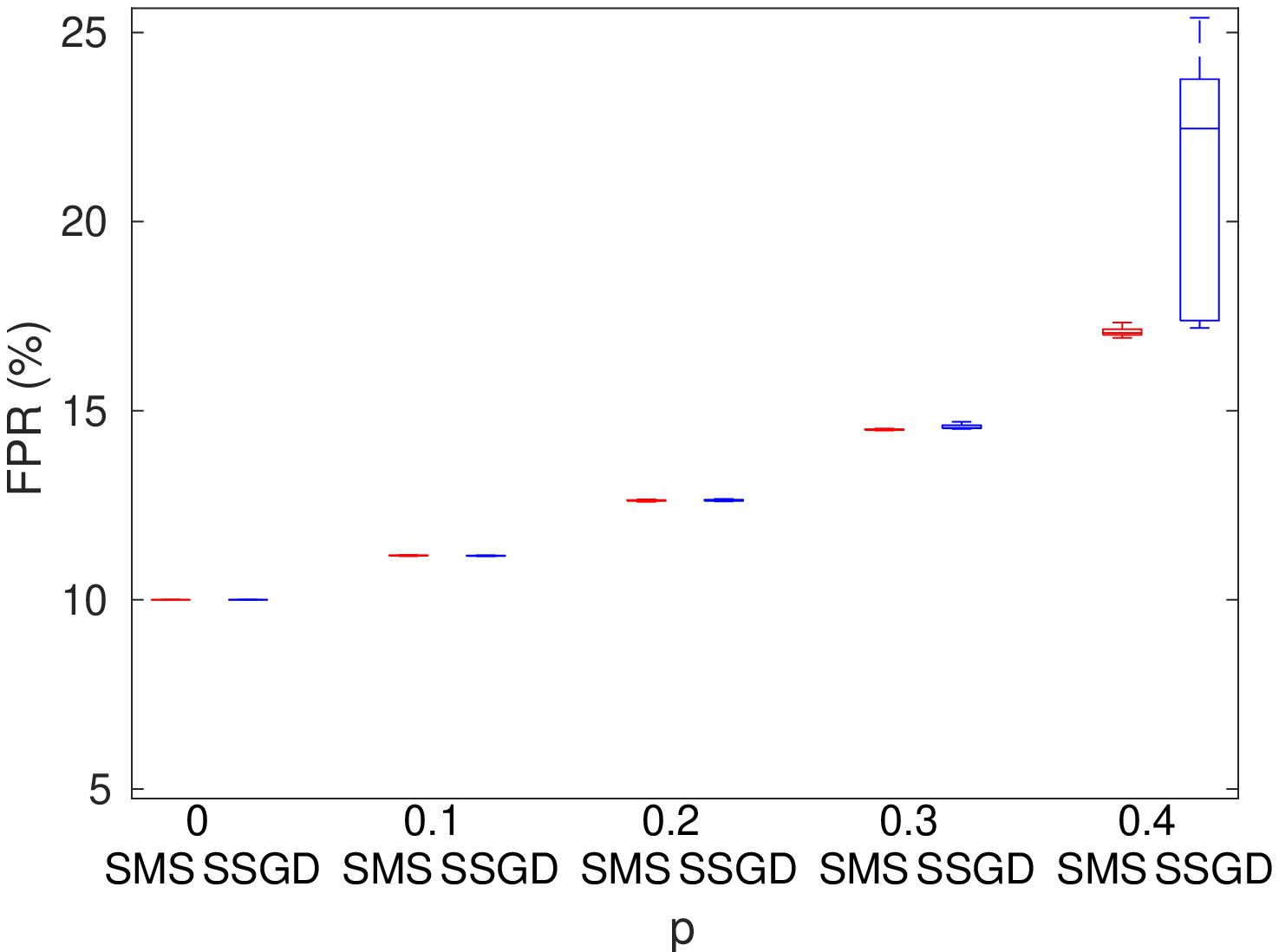}
}
\centerline{\includegraphics[width=.333\textwidth]{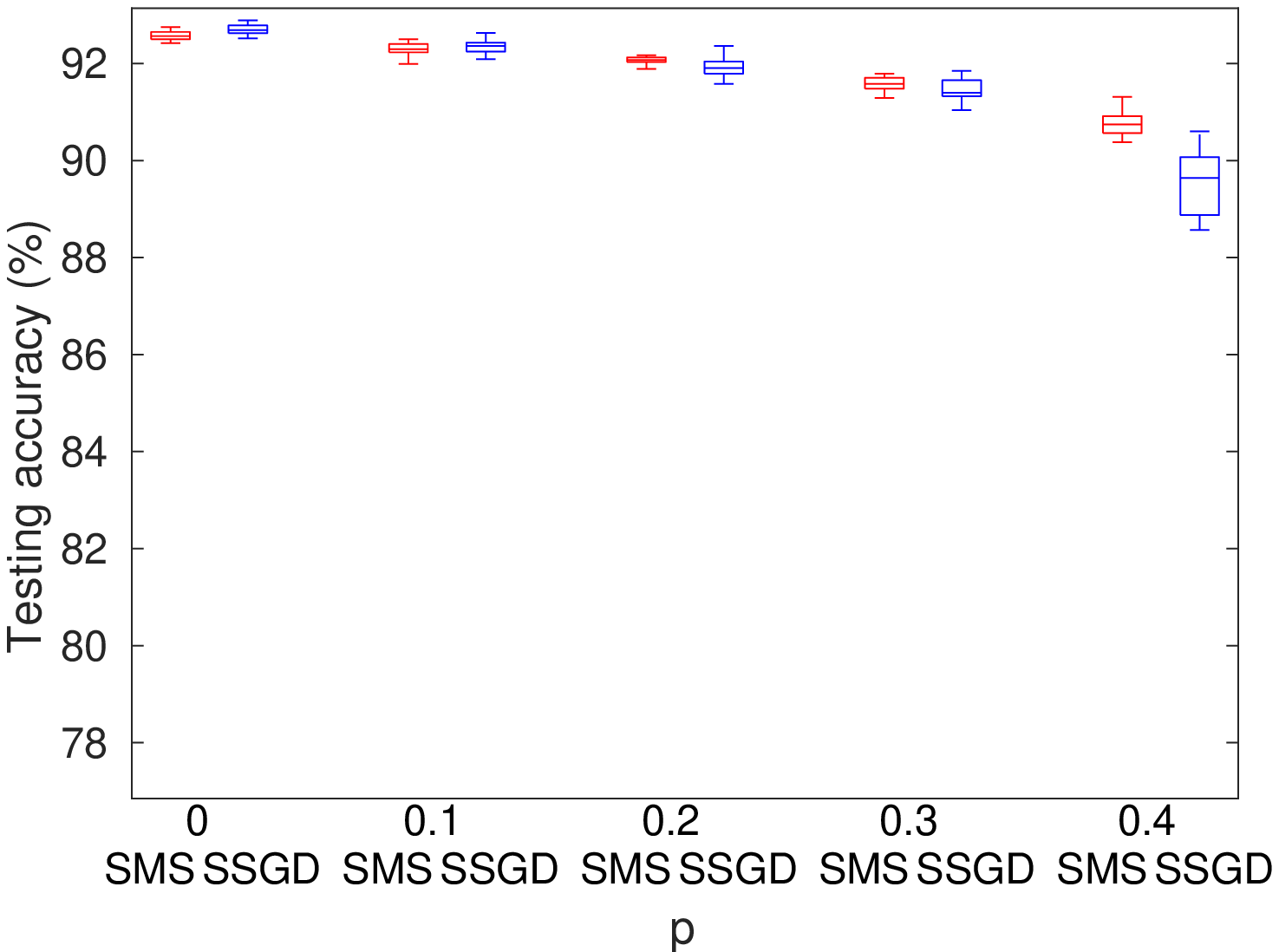}
\includegraphics[width=.333\textwidth]{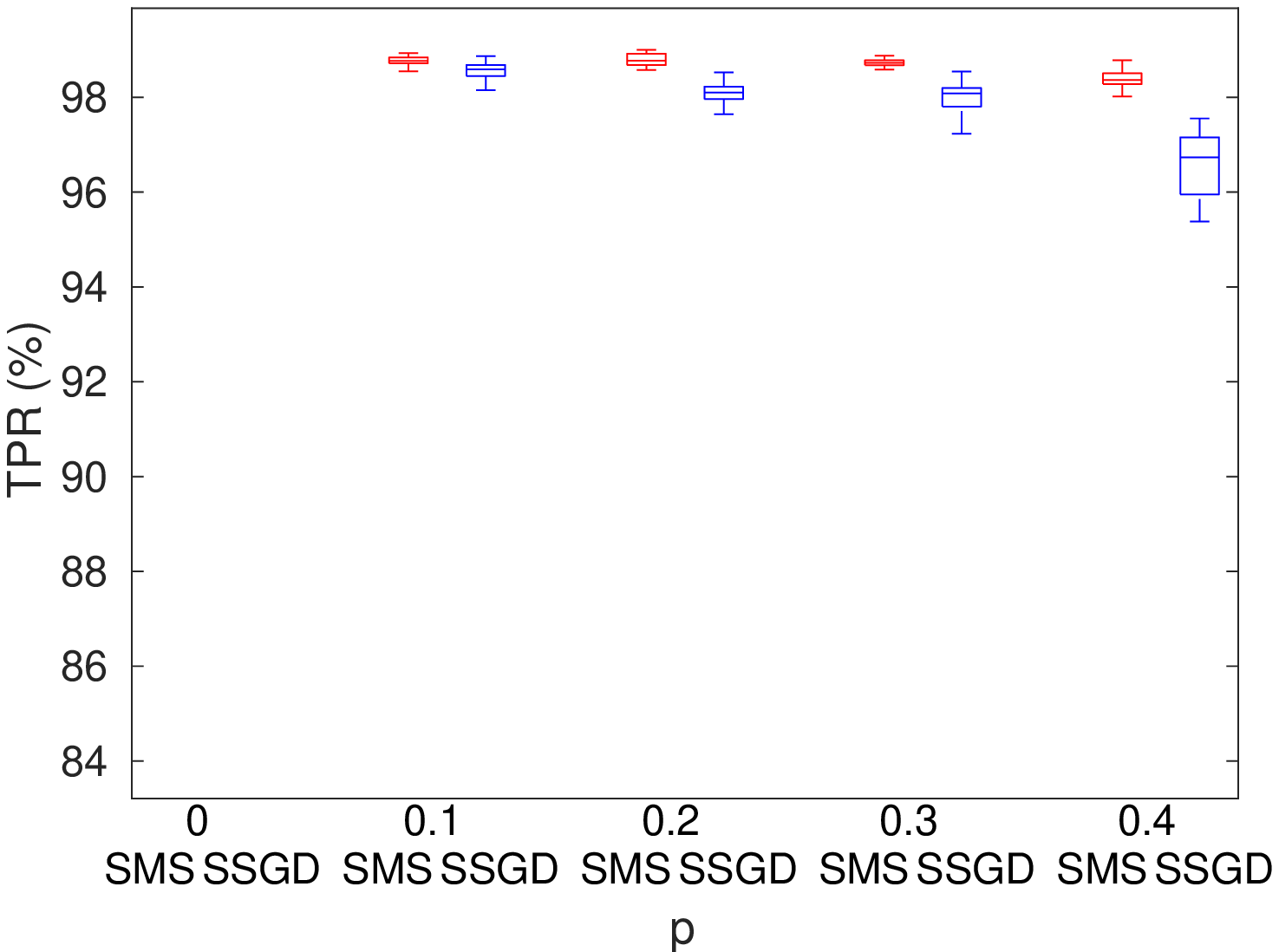}
\includegraphics[width=.333\textwidth]{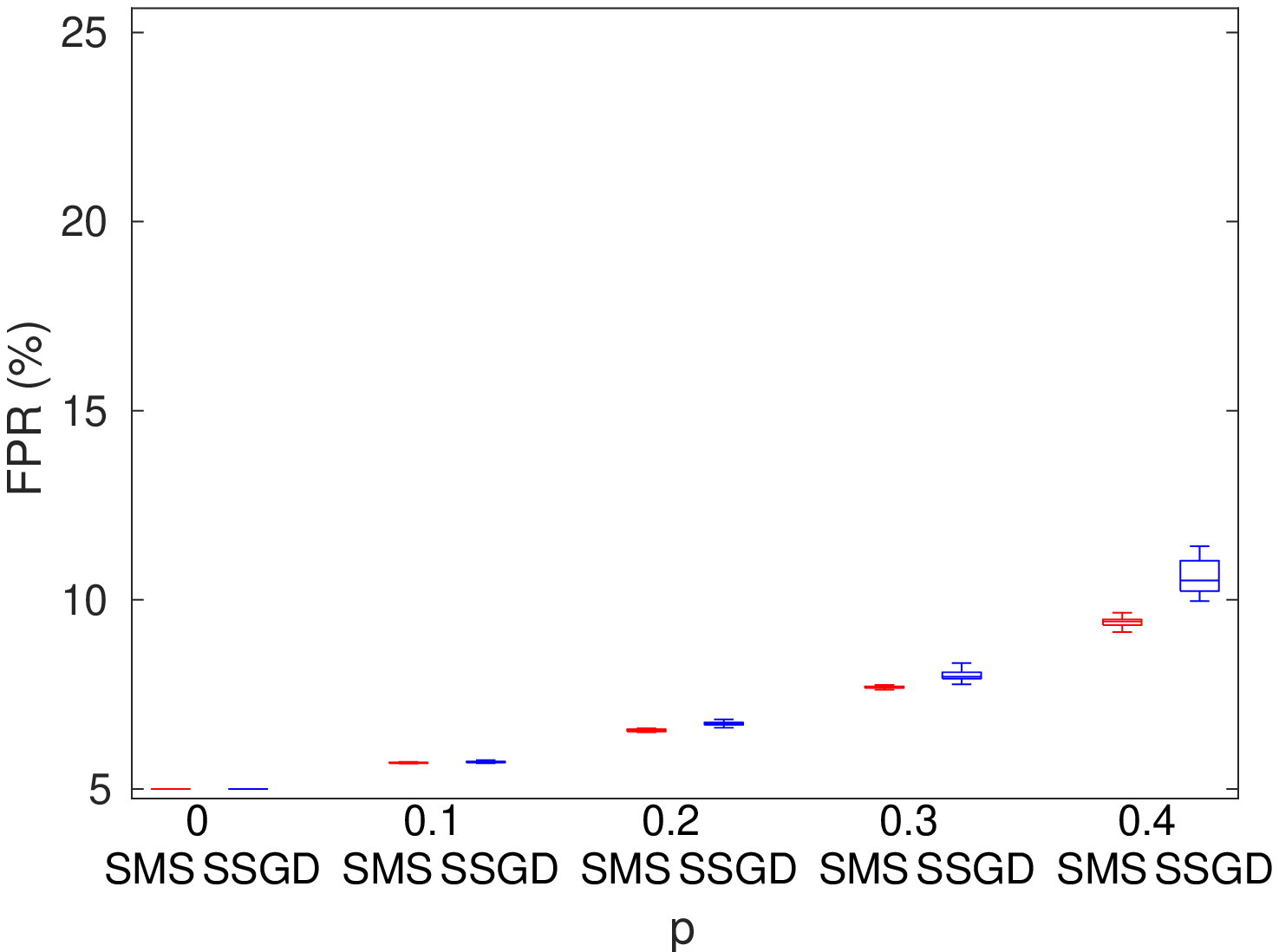}
}
\caption{30 instances of SMS and SSGD on the MNIST classifier 
problem.
\textbf{Top row:} $q = \lfloor(1-p-0.1)N\rfloor$, \textbf{bottom row:} $q = 
\lfloor(1-p-0.05)N\rfloor$.
\textbf{Left:} testing accuracy, 
\textbf{center:} TPR on training dataset, 
\textbf{right:} FPR on training dataset.}
\label{fig:sms_vs_ssgd}
\end{center}
\end{figure*}

For completeness, we show in Figure~\ref{fig:sms_vs_ssgd} the same plots as in 
Figure~\ref{fig:mnist_results}, but we simply replace SGD with the well-tuned 
SSGD (i.e., with $\alpha=0.1$ and $s=40$). 
Interestingly, for $p\in\{0,0.1,0.2,0,3\}$, SSGD performs mildly better than 
does SMS in the $q=\lfloor (1-p-0.1)N\rfloor$ setting in test accuracy, while 
SMS performs mildly better in identifying outliers in the training 
data. 
For the highest levels of
contamination (i.e., $p=0.4$), however, SSGD achieves an undesirably 
wide range of values in all metrics.
In the $q=\lfloor (1-p-0.05)N\rfloor$ setting, SMS consistently outperforms SSGD in all metrics. 

Although we do not illustrate this result here, we remark that a well-tuned 
version of SSGD can perform as well as SMS on the previously discussed 
regression problems,
both in terms of solution time and identification of outliers.

%

\end{document}